\documentclass[11pt,reqno]{amsart}
\theoremstyle{definition}
\newtheorem{rem}{Remark}[section]
\newtheorem{exm}{Example}[section]
\newcommand{\be}{\begin{equation}}
\newcommand{\ee}{\end{equation}}
\newcommand{\bea}{\begin{eqnarray}}
\newcommand{\eea}{\end{eqnarray}}
\newcommand{\beb}{\begin{eqnarray*}}
\newcommand{\eeb}{\end{eqnarray*}}
\usepackage{amssymb,url,amsfonts,amsthm,setspace,indentfirst,mathrsfs}
\usepackage{latexsym,amssymb,amsthm,amsmath}
\usepackage{graphicx}
\usepackage{minibox}
\usepackage{enumitem}
\numberwithin{equation}{section}
\begin{document}
%
\title[On semisymmetric and pseudosymmetric type manifolds]{ On the existence of various generalizations of semisymmetric and pseudosymmetric type manifolds}
\author[A. A. Shaikh]{Absos Ali Shaikh}
\date{\today}
\address{\noindent\newline$^1$ Department of Mathematics,
	\newline University of Burdwan, 
	\newline Golapbag, Burdwan-713104,
	\newline West Bengal, India} 
\email{aask2003@yahoo.co.in, aashaikh@math.buruniv.ac.in}

\dedicatory{}
\begin{abstract}
The objective of the paper is to investigate a sequential study of different generalizations of semisymmetric and pseudosymmetric manifolds with their proper existence by several spacetimes. In the literature of differential geometry, there are many generalizations of such notions in various directions by involving different curvature tensors.  In this paper, we have systematically and consecutively reviewed various generalized notions of semisymmetry, such as, Ricci semisymmetry, conformal semisymmetry, pseudosymmetry, Ricci pseudosymmetry, Ricci generalized pseudosymmetry, conformal pseudosymmetry, Ricci generalized Weyl pseudosymmetry and also many other semisymmetry type conditions. Most importantly, we have exhibited a plenty of suitable examples to examine the proper existence of such geometric  structures, and they are physically significant as several spacetimes admit such geometric structures. By considering an immersion of Schwarzschild black hole metric, we have also provided a new example of Ricci pseudosymmetric manifolds. Finally, the interesting character of Weyl projective curvature tensors of type $(1,3)$ and $(0,4)$ are exhibited with their interesting examples.

\end{abstract}
%
\maketitle
\section{Introduction}
In differential geometry, the notion of symmetry carries out a significant role as it describes the shape and detailed geometric information of a semi-Riemannian manifold. In 1926, Cartan \cite{Cart26} first introduced such notion locally, defined by $\nabla R=0$, which was a generalization of the manifold of constant curvature. We mention that the Otsuki spacetime \cite{Otsuki00, Otsuki01} is  of constant curvature and consequently locally symmetric, but the locally symmetric manifolds, such as, AdS spacetime \cite{AM84, SKMHH03}, Nariai spacetime \cite{Nariai50, Nariai51, SAAC20N} are not  of constant curvature. Again, by using second order covariant differentiation, Cartan \cite{Cartan46} defined the idea of semisymmetry by the condition $R\cdot R=0$, whose classification was obtained by Szab\'o \cite{Szabo82, Szabo84, Szabo85} for the Riemannian manifolds only. Naturally, the idea of semisymmetry generalizes the notion of local symmetry. Again, the semisymmetric manifolds form a subclass of the class of Ricci semisymmetric manifolds, which arises due to the weaker curvature condition  $R\cdot S=0$. For detailed study, we refer the readers to go through the articles \cite{ ACDE98, DGHS11, SKgrt}. Moreover, several generalized notions of semisymmetry were established with different semisymmetric type curvature conditions  by appointing other curvature tensors. Again, the idea of recurrent manifold was first introduced by Ruse \cite{Ruse46, Ruse49a, Ruse49b} as an extension of local symmetry. Later Walker \cite{Walker50} investigated such notion rigorously, and during the last eight decades several generalizations of recurrent manifolds have been introduced by debilitating the curvature condition of recurrency. We must suggest the readers to \cite{Chak87,LR89, SP10, SR10, SR11, SRK16,SRK15,SKP03,TB89, TB93,SKH_2011,SK12,SDHJK,MS12,MS13,MS14} for detailed investigation of such generalized notions of semisymmetry. Recently, Shaikh and Chakraborty \cite{SC_review_recurrent} comprehensively reviewed the first order recurrent like structures with their proper examples in several spacetimes.\\

\indent During the investigation of geodesic mappings as well as totally umbilical submanifolds of semisymmetric manifolds, Adam\'ow and Deszcz \cite{AD83} conceived the idea of pseudosymmetric manifolds, which is also a generalized notion of semisymmetric manifolds. A manifold is \textit{pseudosymmetric} if the tensors $R\cdot R$ and $Q(g,R)$ are linearly dependent i.e., $$R\cdot R=f_RQ(g,R)$$ holds on $\mathscr U_R=\left\lbrace x\in M :  R-\frac{\kappa}{2n(n-1)}(g\wedge g)\ne 0 \ \text{at} \ x\right\rbrace $ where $f_R$ is some scalar function on the preceding set. An extension of pseudosymmetric manifolds form Ricci pseudosymmetric manifolds. A manifold is called  \textit{Ricci pseudosymmetric} if the tensors $R\cdot S$ and $Q(g,S)$ are linearly dependent or equivalently $$R\cdot S=f_SQ(g,S)$$ holds on $\mathscr U_S=\left\lbrace x\in M : S-\frac{\kappa}{n} g\ne 0 \ \text{at}\ x \right\rbrace $ where $f_S$ is some function on $\mathscr U_S$. We mention that this type of symmetry was independently studied by Mike\v{s} and his co-authors in \cite{MIKES76, MS94, MIKES88, MIKES92, MIKES96, MSV15, MVH09} and they named such notion as generalized semisymmetry. Every pseudosymmetric manifold is Ricci pseudosymmetric but not conversely. However, both the notions are equivalent for every $3$-dimensional manifold as well as for every manifolds of which Weyl conformal curvature vanishes identically. We prescribe the reader to go through \cite{Deszcz89, DGHS11,  DHS99, SKgrt} for the study of equivalency of these notions. Further, the imposition of other curvature tensors in the definition of pseudosymmetry leads to various pseudosymmetric type manifolds, which generalizes pseudosymmetry in different directions. \\
\indent The goal of this paper is to discuss various  generalizations of semisymmetric manifolds and explore their proper existence. Thus, this paper may be treated as a review article on the geometric structures of semisymmetry as well as pseudosymmetry. Several ways of generalization of semisymmetric manifolds are presented in the following figures (see, Figure \ref{Fig1} and Figure \ref{Fig2}). In the following figures the rectangular boxes are used to represent the different classes of semisymmtric and pseudosymmetric manifolds, and the ways of generalizations are shown by implication symbols.

\begin{center}
\begin{figure}[h]%
\centering
\includegraphics[width=0.7\textwidth]{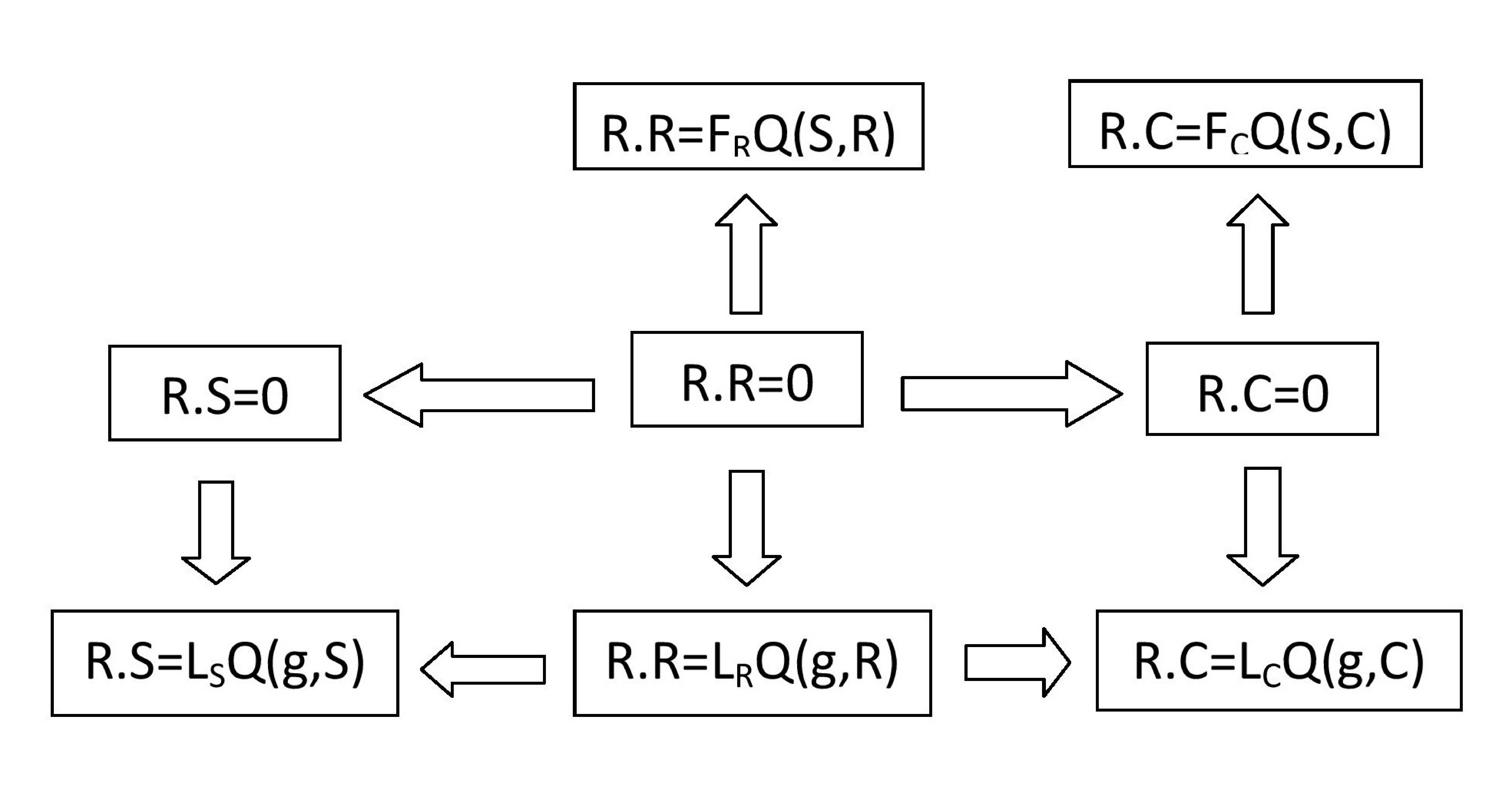}
\caption{Some generalizations of semisymmetric manifold.}\label{Fig1}
\end{figure}
\end{center}

%
\begin{center}
\begin{figure}[h]%
\centering
\includegraphics[width=0.7\textwidth]{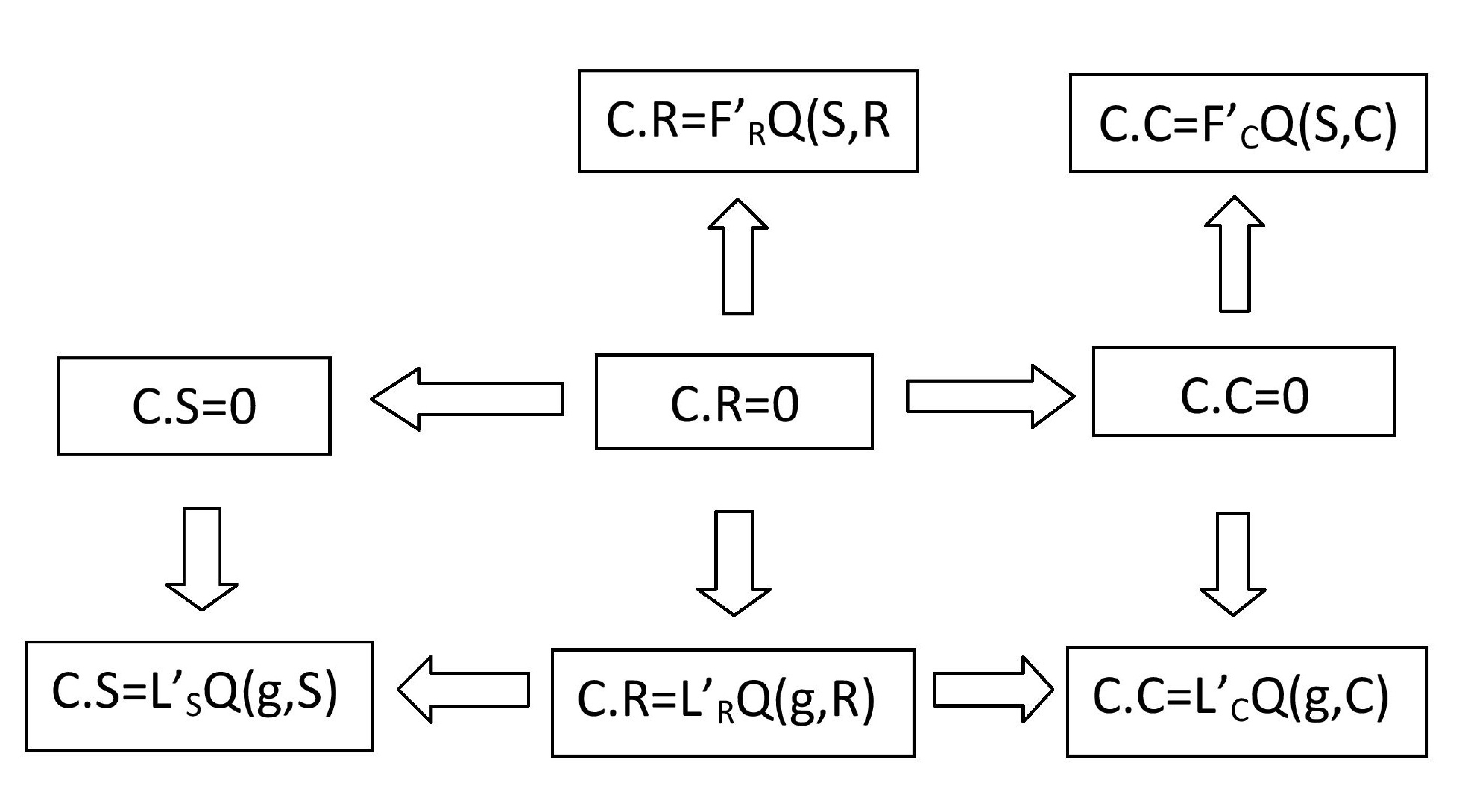}
\caption{Some generalizations of semisymmetric manifold due to Weyl conformal curvature.}\label{Fig2}
\end{figure}
\end{center}
%
 \indent In Figure \ref{Fig1} four different ways of generalization of semisymmetric manifold ($R\cdot R=0$) are shown and these are (i) Ricci semisymmetric manifold ($R\cdot S=0$), (ii) conformally semisymmetric manifold ($R\cdot C=0$), (iii) pseudosymmetric manifold ($R\cdot R=L_RQ(g,R)$) and (iv) Ricci generalized pseudosymmetric manifold ($R\cdot R=F_RQ(S,R)$). Again, Ricci semisymmetric manifold is further extended to Ricci pseudosymmetric manifold ($R\cdot S=L_SQ(g,S)$), and conformally semisymmetric manifold is generalized to (i) conformally pseudosymmetric manifold ($R\cdot C=L_CQ(g,C)$) and (ii) Ricci generalized Weyl pseudosymmetric manifold ($R\cdot C=F_CQ(S,C)$). Both the classes of Ricci pseudosymmetry and conformal pseudosymmetry include the class of pseudosymmetric manifolds. Similar to Figure \ref{Fig1}, the ways of generalization of semisymmetry due to Weyl conformal curvature ($C\cdot R=0$) are shown in Figure \ref{Fig2}.\\
 \indent In the sequel,  we examine the proper existence of the classes of  manifolds given in the following figure (see, Figure \ref{Fig3}). In fact, several spacetimes from the literature of general relativity and cosmology, are found to be the evidences for the existences of such geometric structures.
\begin{center}
\begin{figure}[h]%
\centering
\includegraphics[width=0.7\textwidth]{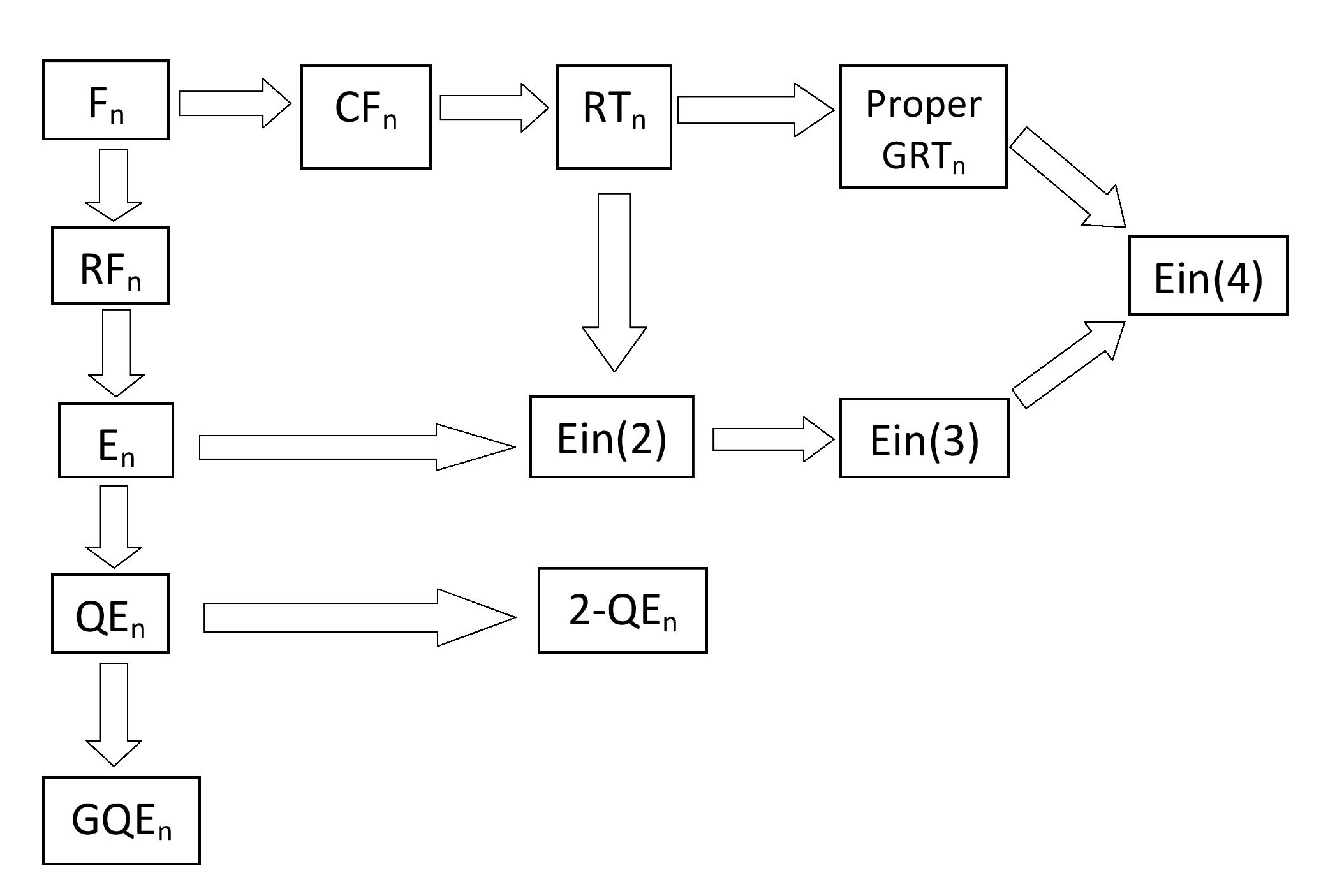}
\caption{Generalizations of flat manifold in different direction.}\label{Fig3} 
\end{figure}
\end{center}
\indent In Figure \ref{Fig3}, two different ways of generalizations of flat manifold ($F_n$) are shown and these are (i) conformally flat manifold ($CF_n$) and (ii) Ricci flat manifold ($RF_n$). Then conformally flat manifold is extended to the Roter type manifold ($RT_n$), the class of proper generalized Roter type manifold ($GRT_n$) includes entirely the Roter type manifolds, and Einstein manifolds of level $4$ ($Ein(4)$) is a further extension of generalized Roter type manifolds. In sequential study, Ricci flat manifold is generalized  to Einstein manifold ($E_n$), Einstein manifold is extended to quasi-Einstein manifold ($QE_n$) and quasi-Einstein manifold is further generalized to generalized quasi-Einstein manifold ($GQE_n$) in the sense of Chaki \cite{Chaki01}. Again, the notion of pseudo quasi-Einstein manifold ($PQE_n$) by Shaikh \cite{Shaikh_2009_PQE} is another generalization of quasi-Einstein manifold.

Einstein manifold is also generalized to Einstein manifold of levels $2$, $3$ and then level $4$. The class of $2$-quasi-Einstein manifolds is again another worthy generalization of quasi-Einstein manifold. \\
%
\indent The outline of this paper is as follows: section $2$ deals with some introductory facts of various tensors. In section $3$, we present the examples which ensure the properness of the classes of manifolds given in Figure \ref{Fig3}. In section $4$, we discuss  various semisymmetric and pseudosymmetric manifolds and examine their proper existence by suitable examples. Section 5 is concerned with the existence of various Weyl symmetric and pseudosymmetric manifolds due to Weyl conformal curvature tensor. In section $6$, we exhibit  various spacetimes satisfying several kinds of pseudosymmetric curvature conditions. The last section deals with the exhibition of the eccentric characteristics of Weyl projective curvature tensors of type (1,3) and (0,4) for the semisymmetric as well as  pseudosymmetric manifolds, and verified the distinction of such tensors by suitable examples. 

\section{Preliminaries}
\indent In this section we define some tensors both in global and local forms on an $n$-dimensional smooth and connected semi-Riemannian manifold $M$. Let $\nabla$ be the Levi-Civita connection on $M$ and $g$ be a semi-Riemannian metric on $M$. We consider the following spaces:
\beb
C^\infty(M)\ &:&\ \mbox{the ring of smooth functions on $M$,}\\
\chi(M)\ &:&\ \mbox{the space of all smooth vector fields on $M$,}\\
\Xi(M)\ &:& \ \mbox{the space of all endomorphisms on} \ \chi(M),\\
\mathcal T^r_k(M) \ &:& \ \mbox{the space of all smooth tensor fields of type}\ (r,k).
\eeb
\indent For $E, D \in \mathcal T^0_2(M)$, the tensor $E\wedge D \in \mathcal T^0_4(M)$ is the Kulkarni-Nomizu product of $E,D$, and it is defined as (\cite{DGHS11}, \cite{Glog02}, \cite{SKA16}, \cite{SRK16})
\beb
(E \wedge D)(\zeta_1,\zeta_2,\zeta_3,\zeta_4)&=&E(\zeta_1,\zeta_4)D(\zeta_2,\zeta_3) + E(\zeta_2,\zeta_3)D(\zeta_1,\zeta_4)\\\nonumber
&-&E(\zeta_1,\zeta_3)D(\zeta_2,\zeta_4) - E(\zeta_2,\zeta_4)D(\zeta_1,\zeta_3),
\eeb
where $\zeta_1,\zeta_2,\zeta_3,\zeta_4 \in \chi(M)$. Locally $E\wedge D$ is written as 
$$(D\wedge E)_{abcd}=D_{ad}E_{bc} + D_{bc}E_{ad} -D_{ac}E_{bd} - D_{bd}E_{ac}.$$
\indent Again, for a symmetric tensor $N \in \mathcal T^0_2(M)$, we get $(\zeta_1\wedge_N \zeta_2) \in \Xi(M)$, $\mathcal{N} \in \Xi(M)$ (\cite{DD91}, see also \cite{SKgrt}, \cite{SKgrtw}) which is given in below:
$$(\zeta_1\wedge_N \zeta_2)\zeta = N(\zeta_2,\zeta)\zeta_1-N(\zeta_1,\zeta)\zeta_2, \ \ g(\mathcal N \zeta_1,\zeta_2) = N(\zeta_1,\zeta_2).$$
For $L\in \mathcal T^0_4(M)$ and $\zeta,\eta \in \chi(M)$, the associated tensor $\mathcal L\in \mathcal T^1_3(M)$ and the associated curvature operator $\mathscr L(\zeta,\eta) \in \Xi(M)$ are respectively given by
$$g(\mathcal{L}(\zeta,\eta)\zeta_1,\zeta_2)=L(\zeta,\eta,\zeta_1, \zeta_2) \ \mbox{and}$$
$$\mathscr{L}(\zeta,\eta)(\zeta_1)=\mathcal L(\zeta,\eta)\zeta_1.$$

\indent  Let $H\in \mathcal T^0_4(M)$ be such that 
\beb
\mathop{\mathcal S}_{\zeta_1,\zeta_2,\zeta_3} H(\zeta_1,\zeta_2,\zeta_3,\zeta_4)=0,
\eeb
$$H(\zeta_1,\zeta_2,\zeta_3,\zeta_4) + H(\zeta_2,\zeta_1,\zeta_3,\zeta_4)=0 \ \ \mbox{and}$$
$$H(\zeta_1,\zeta_2,\zeta_3,\zeta_4)=H(\zeta_3,\zeta_4,\zeta_1,\zeta_2),$$
where $\mathcal S$ stands for cyclic sum over the vector fields $\zeta_1, \zeta_2, \zeta_3$. Then $H$ is a generalized curvature tensor (\cite{DGHS11}, \cite{SDHJK15}, \cite{SK14}) and it is mentioned that $E\wedge D$ is a generalized curvature tensor if $E,D \in \mathcal T^0_2(M)$ both are symmetric.
Again a generalized curvature tensor is said to be proper if 
\beb
\mathop{\mathcal S}_{\zeta_1, \zeta_2, \zeta_3} (\nabla_{\zeta_1}H)(\zeta_2, \zeta_3, \zeta_4, \zeta)=0.
\eeb
Now, we consider the tensor $H\in \mathcal{T}^0_4(M)$ written in terms of the linear combination of $R$, $g\wedge g$ and $g\wedge S$ as :
\bea\label{LT}
H=a_1 R + a_2(g\wedge g) + a_3 (g\wedge S).
\eea
From \eqref{LT} we get the following tensors:
\beb
&&\mbox{the Gaussian curvature} \ G \  \mbox{for} \ a_1=0,\ a_2=\frac{1}{2}  ,\  a_3=0, \\
&&\mbox{the conformal curvature} \ C \  \mbox{for} \ a_1=1,\ a_2=\frac{\kappa}{2(n-1)(n-2)},\  a_3=-\frac{1}{n-2},\\
&&\mbox{the conharmonic curvature} \ K \  \mbox{for} \ a_1=1,\ a_2=0,\  a_3=-\frac{1}{n-2},\\
&&\mbox{the concircular curvature} \ W \  \mbox{for} \ a_1=1,\ a_2=-\frac{\kappa}{2n(n-1)},\  a_3=0.
\eeb
Now the local forms of the Riemann curvature tensor, the Ricci curvature and the scalar curvature are given by
$$R_{abcd} = g_{a\alpha}(\partial_{d}\Gamma^{\alpha}_{bc} - \partial_c \Gamma^{\alpha}_{bd} + \Gamma^{\beta}_{bc}\Gamma^{\alpha}_{\beta d} - \Gamma^{\beta}_{bd}\Gamma^{\alpha}_{\beta c}),\ S_{ab}=G^{\alpha \beta}R_{\alpha ab \beta}\ \ \mbox{and} \ \ \kappa=g^{ab}S_{ab},$$
and by these components one can easily find the local forms of all the curvatures mentioned above. Here $\Gamma^{\alpha}_{ab}$ are $2$nd kind Christoffel symbols for the Levi-Civita connection $\nabla$ and $\partial_\alpha=\frac{\partial}{\partial x^\alpha}$. We also deal with another curvature $P \in \mathcal T^0_4(M)$, the projective curvature, given by:
$$P_{abcd}=R_{abcd} - \frac{1}{n-1}\left( g_{ad}S_{bc} - g_{bd}S_{ac}\right).$$ 
For detailed investigation of the 2nd order curvature restricted geometric structure, we refer the reader to go through the work of Shaikh and Kundu \cite{Shaikh_Kundu_IJGMMP_2018}.

Let $U \in \mathcal T^0_k(M)$, $k\ge 1$. Then we define the tensor $(L\cdot U)\in \mathcal T^0_{k+2}(M)$ by  (see, \cite{DGHS98, DH03, SK14})
\beb
(L\cdot U)_{i_1i_2\cdots i_kab}=-g^{\alpha \beta}\left[L_{abi_1\alpha}U_{\beta i_2\cdots i_k}+\cdots + L_{abi_k \alpha}U_{i_1i_2\cdots \beta} \right].
\eeb
Also the tensor $Q(N,U)\in \mathcal T^0_{k+2}(M)$, called Tachibana tensor \cite{DGPSS11, SDHJK15, Tachibana74}, is defined as:
\beb
Q(N, U)_{i_1 i_2\cdots i_k a b} &=& N_{bi_1}U_{a i_2\cdots i_k} + \cdots + N_{bi_k}U_{i_1 i_2\cdots a}\\
& & -N_{ai_1}U_{b i_2\cdots i_k} - \cdots - N_{ai_k}U_{i_1 i_2\cdots b}.
\eeb
For $L=R$ and $U=R$, $C$, $K$, $W$, $P$ and $S$ we get the $(0,6)$ type tensors $R\cdot R$, $R\cdot C$, $R\cdot K$, $R\cdot W$, $R\cdot P$ and $R\cdot S$ respectively. Replacing $L$ by $C$, $K$, $W$, $P$ we can get the corresponding $(0,6)$ type curvature tensors. Again, for $N=g$ and $U=R$, $C$, $K$, $W$, $P$ and $S$ we get the Tachibana tensors $Q(g, R)$, $Q(g, C)$, $Q(g,K)$, $Q(g,W)$, $Q(g,P)$ and $Q(g,S)$ respectively. Again for $N=S$ we get the corresponding Tachibana tensors.

\section{Proper existence of the classes of manifolds presented in Figure \ref{Fig3}}
%
\subsection{Flat, Ricci flat and conformally flat manifold}
An $n$-dimensional semi-Riemannian manifold, which is locally isomorphic to an Euclidean manifold, is called a flat manifold. In terms of curvature a semi-Riemannian manifold is flat if and only if $R=0$ i.e., its Riemann curvature vanishes identically.
\begin{exm}\label{ExmF_n}
	The Minkowski spacetime \cite{SKMHH03} given by the metric
	\beb
	ds^2=-dt^2+dx^2+dy^2+dz^2
	\eeb
	is a flat manifold.
\end{exm}
  Again a manifold with $S=0$ is called a Ricci flat manifold. In general relativity, a vacuum solution of Einstein equations is a Ricci flat manifold. Evidently, every flat manifold is Ricci flat but not conversely, and it is given in the following example. 
  \begin{exm}\label{ExmRF_n}
     Kasner spacetime \cite{Kasner21} is a vacuum solution of general relativity and its metric is given by
  	\beb   ds^2=dt^2+t^{2a_1}dx^2+t^{2a_2}dy^2+t^{2a_3}dz^2
  	\eeb 
  	where $a_1^2+a_2^2+a^2_3=1$ and          $a_1+a_2+a_3=1$. It is a Ricci          flat manifold (i.e., $S=0$), but it is not           flat as $R\ne 0$ and the non-zero components of $R$ are given as follows:
   $$R_{1k1k}=-a_{k-1}(a_{k-1}-1)t^{2(a_{k-1}-1)},\ \ \text{ for } k=2,3,4$$
   $$R_{ijij}=-a_{i-1}a_{j-1}t^{2(a_{i-1}+a_{j-1}-1)},\ \ \text{ for } (i,j)=(2,3),(2,4),(3,4).$$
   \end{exm}

   \begin{exm}
     In terms of the outgoing Eddington–Finkelstein coordinates $(u,r,\theta,\phi)$, the Schwarzschild black hole spacetime  for the exterior field is given by \cite{GP09}
  	\beb
  	ds^2=-\left( 1-\frac{2m}{r}\right)       du^2 -2dudr + r^2\left(d\theta^2       + \sin^2\theta d\phi^2 \right), 
  	\eeb
      where $m$ is a constant. Then the spacetime is Ricci flat as $S=0$, but it is not  flat  as the non-zero components of $R$ are given by:
       $$R_{1212}=-\frac{2m}{r^3}=-\frac{2}{r^2}R_{1323}=-\frac{2}{r^2 \sin^2\theta}R_{1424},$$
       $$R_{1313}=\frac{m(r-2m)}{r^2}=\frac{1}{\sin^2\theta}R_{1414},\ \ R_{3434}=2mr\sin^2\theta.$$
\end{exm}
   An extension of flat manifold is conformally flat manifold ($C=0$) \cite{SC67}. This class perceives great importance in general relativity and cosmology as Weyl conformal curvature is relativistic equivalent to tidal forces of spacetimes. Naturally every flat manifold is conformally flat, but not conversely.
   \begin{exm}\label{ExmCF_n}
   	The warped product $M_1\times_f M_2$ of an $1$-dimensional manifold $(M_1,\bar g)$, $\bar g_{11}=-1$, with warping function $f$ and a $3$-dimensional Riemannian space $(M_2,\tilde g)$ of constant curvature is called a Robertson-Walker spacetime \cite{ARS95, Neill83,S98, SKMHH03}. The Weyl conformal curvature for this spacetime vanishes identically.
   	\end{exm}

A manifold is of constant curvature if its curvature tensor $R$ can be expressed as $R=k(g \wedge g)$. Every flat manifold is obviously of constant curvature, but not conversely.

\begin{exm}
	The Otsuki spacetime \cite{Otsuki00, Otsuki01}, in terms of $(x_1,x_2,x_3,x_4)$ coordinates, is given by
	$$ds^2=\frac{1}{dx_4^2}\left\{-\frac{1}{1+ax_4^2}dx_4^2 + \sum_{i,j=1}^{3}\left(\delta_{ij}-\frac{ax_i x_j}{1+ar^2}\right) dx_i dx_j     \right\},$$
	where $a=constant$ and $r^2=\sum_{k=1}^{3}x_k^2$, and it is a manifold of constant curvature, which is not flat.
\end{exm}

%
%
\subsection{Roter type and generalized Roter type manifold}
%
\indent Generalizing the conformally flat manifold, Deszcz \cite{Deszcz03} presented the Riemann curvature in terms of some lower order tensors as
\be\label{RT}
R=\mu_1(g\wedge g) + \mu_2 (g\wedge S) + \mu_3 (S\wedge S). \notag
\ee
This condition is known as Roter type ($RT_n$) condition. Very recently, in \cite{DDH21} this condition is proved to be equivalent with pseudosymmetry condition on spacetimes. It seems that Reissner-Nordstr\"om metric \cite{Kowa06} first appeared with Roter type spacetime. We present another spacetimes satisfying such condition.
\begin{exm}\label{ExmRT_n}
 Rabinson-Trautman spacetime \cite{RT62} represents an isolated gravitationally radiating system and the metric is given by
    \beb
    ds^2=-2\left(U^0 - 2\gamma^0r - \Psi^0_2r^{-1}\right)du^2 + 2dudr - \frac{r^2}{2P^2} d\xi d\bar{\xi},
    \eeb 
    where $U^0$, $\gamma^0$,$\Psi^0_2$ are constants and $P$ is a nowhere vanishing function of $\xi$ and $\bar \xi$. Also, the metric can be written as
    \begin{equation}
    ds^2=-2\left(a-2br-\frac{q}{r}\right)dt^2+2dtdr-\frac{r^2}{f^2}(dx^2+dy^2),
    \end{equation}
    where $a,b,q$ are constants and $f$ is a nowhere vanishing function of $x,y$. Recently, Shaikh et al. \cite{SAA18} proved that the  Rabinson-Trautman spacetime admits Roter type condition as given in below:
    $$R=a_{11}g\wedge g + a_{12}g\wedge S + a_{22}S\wedge S,$$
    where $a_{ij}$'s are given by
    \begin{eqnarray}
        a_{11}&=&-\frac{1}{r^3}\left[q+\frac{4br(6aq+8b^2r^3-24bqr-3qF)}{((-2a + 4br + F)^2}      \right],\nonumber\\
        a_{12}&=& \frac{4abr^2 + 6aq + 8b^2r^3 - F(2br^2 + 3q)-36bqr}{r(-2a+4br+F)^2},\nonumber\\
        a_{22}&=&\frac{r(2ar-6q-Fr)}{2(-2a+4br+F)^2},\nonumber
    \end{eqnarray}
    with $F=f_{xx}^2+f_{yy}^2-f(f_{xx}+f_{yy})$ provided $ F-2a+ 4br \neq 0$ everywhere.
\end{exm}
A semi-Riemannian manifold  satisfying the condition
\be
R=J_1(g\wedge g) + J_2 (g\wedge S) + J_3(g\wedge S^2) + J_4(S\wedge S) + J_5(S\wedge S^2) + J_6(S^2\wedge S^2), \notag
\ee
is called a generalized Roter type ($GRT_n$) manifold. A generalized Roter type manifold, which is not $RT_n$, is called a  proper generalized Roter type manifold. Such notion was investigated in \cite{DGJPZ13, DHJKS14, Sawicz06}, but the term ``generalized Roter type'' was first used in \cite{SDHJK15}.

\begin{exm}\label{ExmGRT_n}
         Kantowski-Sachs spacetime \cite{KS66} is a homogeneous and anisotropic model of the universe and its metric is given by 
        \be
        ds^2=dt^2-X^2(t)dr^2-Y^2(t)[d\theta^2+\sin^2(\theta) d\phi^2].\notag
        \ee
        In \cite{SC20} it is proven that this spacetime fulfills proper $GRT_4$ condition.
\end{exm}

\begin{exm}\label{exm-melvin}
    The Melvin magnetic metric is static cylindrically symmetric solution of Einstein-Maxwell equation, representing a bundle of parallel magnetic lines of force, under which  their mutual gravitational attraction remain in equilibrium. The Melvin magnetic metric is given by
    \beb
    ds^2=\left(1 + \frac{B_0r^2}{4}\right)(-dt^2 + dr^2 + dz^2) + \frac{r^2}{\left(1 + \frac{B_0r^2}{4}\right)} d\phi^2,
    \eeb
     $B_0$ being magnetic field at $r=0$ axis. The Weyl form of the Melvin type metric can be written as
     $$ds^2=e^{2f(r)}(-dt^2+dr^2+dz^2)+r^2 e^{-2f(r)}d\phi^2.$$     
     In \cite{SAAC20melvin}, it is shown that such metric fulfills generalized Roter type condition given as follows:
     $$R=a_{12}g\wedge S + a_{13}g\wedge S^2 + a_{22}S\wedge S+a_{23}S\wedge S^2,$$
     where the $a_{22}$'s are given by
    \begin{eqnarray}
        && a_{12}= \frac{f'}{2\lambda}+\frac{r(\lambda e^{2t}+f'')}{2(re^{2t})}, \ \ \ 
        a_{13}=\frac{e^4r^2ff''}{4(r e^2f \mu+f')\lambda^2}, \nonumber \\
        &&a_{22}=a_{13}\lambda e^{-2f} - \frac{r^2 e^{4f} \mu}{2\lambda^2}, \ \ \ 
        a_{23}=2a_{13}r\mu e^{2f}, \nonumber
    \end{eqnarray}
    with $\lambda=f'+rf''$ and $\mu=\frac{e^{-2f}}{r}(f'-rf'^2)$.
\end{exm}

    Recently Shaikh and his coauthors showed that Vaidya-Bonner spacetime \cite{SDC21}, Lema\^itre-Tolman-Bondi spacetime \cite{SAACD22}, Lifshitz spacetime \cite{SSC19} admit the  $GRT_4$ conditions  described as follows:
    
    \begin{exm}
    Vaidya-Bonner spacetime, whose metric is given by \cite{SDC21}
            \begin{eqnarray}
        	ds^2-\left(1-\frac{2m(u)}{r}+\frac{q^2(u)}{r^2}\right)du^2-2dudr+r^2(d\theta^2+\sin^2\theta d\phi^2),\nonumber
            \end{eqnarray}
            admits the following $GRT_4$ relation:
            \begin{eqnarray}
        	R&=&\frac{3r^4}{4q^4}\left( -rm + q^2\right) S\wedge S + \frac{r^8}{2q^6}\left(rm - q^2 \right) S\wedge S^2\notag \\
        	&+& \left( -\frac{rm}{2q^2} + 1\right) g\wedge S + \frac{r^4}{4q^4}\left(-rm + q^2 \right) g\wedge S^2. \notag
            \end{eqnarray}
        \end{exm}

        \begin{exm}
         The Lema\^itre-Tolman-Bondi spacetime  given in Example \ref{LTB_pseudo} satisfies the following $GRT_4$ condition \cite{SAACD22}:
            \bea
            R=a_{22} S \wedge S + a_{23} S\wedge S^2 + a_{12} g\wedge S + a_{13} g\wedge S^2 \notag
            \eea
            where
            \bea
            a_{22} &=& \frac{B^3B'}{D}[-B^5\ddot{B}'\Psi_1 + BB'\Psi_1\Psi_2(4\ddot{B}(\Psi_1 + B'\ddot{B}) - 7\ddot{B}'\Psi_2)  \notag\\
            & & +B^3\ddot{B}'\Psi_1(\Psi_1^2 - 2B'\ddot{B}\frac{\partial \Psi_2}{\partial r}  
             + 7B'\ddot{B}\Psi_2) + B'^2\Psi_2^2(\Psi_1^2 - B'\ddot{B}'\Psi_2) \notag\\
             & &+ B^4\ddot{B}'^2(\Psi_1^2 - 2B'\ddot{B}'\Psi_2) + B^2(-\Psi_1^2((\frac{\partial \Psi_2}{2\partial r})^2 - 10B'(\frac{\partial \Psi_2}{2\partial r})\ddot{B} \notag\\
             & & + B'\ddot{B}^2) + B'\ddot{B}'\Psi_2^2(7\Psi_1^2 + 3B'\ddot{B}'))],\notag\\
             a_{23} &=&\frac{B^4B'^2}{D}[-B'^2\Psi_2^2\Psi_1 + 2BB'^2\ddot{B}'\Psi_2^2 +2B^3B'\ddot{B}'^2\Psi_2 - B^4\ddot{B}'^2\Psi_1 \notag\\
             & & + B^2\Psi_1(\Psi_1^2 + 4B'\dot{B}\dot{B}'\ddot{B}- 4\mu'B'\ddot{B} + 4B'\ddot{B}'\Psi_2)],\notag\\
             a_{12} &=& \frac{B}{D}[B'^2\ddot{B}'\Psi_2^2(B'\Psi_2 + 3B\Psi_1) + B^4\ddot{B}'^3(B\Psi_1 + 3B'\Psi_2) \notag\\
             & & - B^3\ddot{B}'\Psi_1(\Psi_1^2 - 2B'\ddot{B}\frac{\partial \Psi_2}{\partial r} + 6B'\ddot{B}'\Psi_2) + (\Psi_1 + B'\ddot{B})(5\Psi_1\Psi_2\ddot{B}' \notag\\
             & & - 9\Psi_2B'\ddot{B}\ddot{B}' - 8\Psi_1^2\ddot{B}) + 5B'^2\ddot{B}^2\ddot{B}'\Psi_2],\notag\\
            a_{13} &=& \frac{B^3B'^2}{D}[(B'\Psi_2 + B^2\ddot{B}')\Psi_2\ddot{B}' + B\Psi_1(4\ddot{B}(\Psi_1 + B'\ddot{B}) - 3\Psi_2\ddot{B}')],\notag
            \eea
            with $D=2(B'\Psi_2 + B\Psi_1)^2(-\Psi_1 + B\ddot{B}')(B'\Psi_2 + B(-\Psi_1 + B\ddot{B}'))(B'\Psi_2 + B(\Psi_1 + B\ddot{B}')).$
        \end{exm}

        \begin{exm}
         The Lifshitz spacetime, in Poincar\'e coordinates $(t,r,x,y)$, is given by
            $$ds^2=-\frac{dt^2}{r^{2\mu}}+\frac{1}{r^2}(dr^2+dx^2+dy^2),$$
            where $\mu$ is a constant. Shaikh et al. \cite{SSC19} showed that the Lifshitz spacetime is not Roter type, but it satisfies the following $GRT_4$ condition:
            \begin{eqnarray}
               R&=&a_{11}(g\wedge g) + a_{12} (g\wedge S) + a_{13}(g\wedge S^2) \nonumber\\
               &+& a_{22}(S\wedge S) + a_{23}(S\wedge S^2) + a_{33}(S^2\wedge S^2), \nonumber
            \end{eqnarray}
            where  $a_{ij}$'s are given by
            \begin{eqnarray}
                a_{11}&=&-\frac{c(a\mu+6)+13\mu+8}{ab(a\mu+4)} + \frac{c^2[\mu(2a\mu+4)+2]}{a(a\mu+4)} + \frac{c[\mu(a \mu+3)+1]}{a(a\mu+4)},\nonumber\\
                a_{12}&=&\left[\frac{a\mu+4}{bc}-2c-c^2  \right] ,\ \ \ \ \  \ a_{22}=a_{23}=1,\nonumber\\
                a_{13}&=& \frac{1}{a(a\mu+4)} \left[-\frac{3(c\mu+3)}{bc}+3c+bc(a\mu+2)   \right], \nonumber\\
                a_{33}&=&-\frac{1}{ac(c\mu+4)} \left[ (c\mu+3)+\frac{1}{2b}+1 \right],\nonumber  
            \end{eqnarray}
            and $(a,b,c)=(1+\mu,1-\mu,2+\mu)$  with $\mu\neq -2, -1,1$.
    Hence, Lipshitz spacetime is a proper generalized Roter type spacetime.
\end{exm}

%
\subsection{Einstein, quasi-Einstein and generalized quasi-Einstein manifold}
%
Einstein manifolds play crucial role in semi-Riemannian geometry as well as in general relativity. According to Besse \cite{Bess87}, a semi-Riemannian manifold is called Einstein if its Ricci tensor is linearly dependent with the metric tensor i.e., $S=\frac{\kappa}{n}g$. This manifold generalizes both the manifold of constant curvature as well as Ricci flat manifold.
\begin{exm}
	The Kaigorodov  spacetime metric \cite{Kaigorodov63} is given by
	\beb
	ds^2=(dx^4)^2+e^{2x^4/l}[2dx^1dx^3+(dx^2)^2]\pm e^{-x^4/l}(dx^3)^2.
	\eeb
	This spacetime is a non-Ricci flat Einstein manifold \cite{SDKC19}. It is also not a manifold of constant curvature. 
\end{exm}
Again, from the class of semi-Riemannian manifolds a natural subclass of the Einstein manifolds can be obtained by a curvature condition imposed on their Ricci tensors. Also, the class of quasi-Einstein manifolds is another subclass of the class of semi-Riemannian manifolds, which generalizes the notion of Einstein manifold. A manifold is called quasi-Einstein if its Ricci tensor can be decomposed as
\beb
S=pg + \rho (\sigma \otimes \sigma),
\eeb
where $p, \rho$ are smooth functions and $\sigma$ is the associated non-zero 1-form such that $g(\zeta,X)=\sigma(X)$, $\zeta$ being the unit vector field known as generator of the manifold.
\begin{exm}\label{ExmQE_n}
     Robertson-Walker spacetimes \cite{ARS95, Neill83, S98, SKMHH03} mentioned in $Example$ \ref{ExmCF_n} are not Einstein but quasi-Einstein manifolds.
     \end{exm}

    \begin{exm}
     Van Stockum in $1937$ derived a cosmological model of rotating universe whose metric, named as van Stockum dust metric \cite{Van37}, is given by
        \beb
        ds^2=dt^2 - 2 \omega r^2dt d\phi -(1-\omega^2r^2)r^2d\phi^2 - e^{-\omega^2r^2}(dr^2+dz^2).
        \eeb
        The non-zero components of Ricci curvature (upto symmetry) are given as follows:        $$S_{11}=-2\omega^2e^{\omega^2r^2},\ \ S_{12}=2\omega^3r^2e^{\omega^2r^2},\ \ S_{22}=-2\omega^2e^{\omega^2r^2}r^2(1+\omega^2r^2),$$ $$S_{33}=-2\omega^2=S_{44}\ \text{and the scalar curvature}\ \kappa=4\omega^2r^2e^{\omega^2}.$$
        We note that $S=pg + \rho (\sigma \otimes \sigma)$ is satisfied for $p=2\omega^2e^{\omega^2r^2}$, $\rho=-4\omega^2e^{\omega^2r^2}$ and $\sigma=\left\lbrace 1, -\omega r^2,0,0\right\rbrace$ with $\|\sigma\|=1$. Thus the above metric is quasi-Einstein manifold.
\end{exm}
\begin{rem}
	Lorentzian manifolds with quasi-Einstein Ricci tensor appears with great significance in cosmology as every perfect fluid spacetime is quasi-Einstein and conversely \cite{SYH09}. The van Stockum dust solution given in (ii) of $Example$ \ref{ExmQE_n} is a perfect fluid spacetime.
\end{rem}
Again, Chaki \cite{Chaki01} extended the notion of quasi-Einstein manifold to the notion of generalized quasi-Einstein manifold by the following form of the Ricci tensor
\bea \label{GQE}
S=pg + \rho (\sigma\otimes \sigma) + \bar{\rho}(\sigma\otimes \delta + \delta\otimes \sigma)
\eea
where $p, \rho, \bar{\rho}$ are non-zero scalars and $\sigma, \delta$ are non-zeo one forms. It is easy to note that every quasi-Einstein manifold is generalized quasi-Einstein but not conversely. We present the following instance.
\begin{exm}
     The Robinson-Trautman metric given in $Example$ \ref{ExmRT_n} is not a quasi-Einstein but a generalized quasi-Einstein manifold \cite{SAA18} as it satisfies the relation (\ref{GQE})
        for $p=-\frac{F+8br-2a}{r^2}$, $\rho=\frac{4b(F+8br-2a)}{\eta^2r}$, $\bar\rho={F+8br-2a}$, $\sigma=\left\{   \frac{q-ar}{\eta r^3},\frac{1}{\eta r^2}, 0,0 \right\}$ and $\delta=\left\{ \eta,0, 0,0 \right\}$, where $\eta$ is an arbitrary constant and $F=f_{xx}^2+f_{yy}^2-f(f_{xx}+f_{yy})$.
        \end{exm}

        \begin{exm}
         Kantowski-Sachs metric given in $Example$ \ref{ExmGRT_n} is again a non-quasi-Einstein generalized quasi-Einstein manifold \cite{SC20} as it obeys the relation (\ref{GQE}) for
         
         \be
         p=\frac{fX\dot Y^2+fXY\ddot Y+fY\dot X \dot Y-X\ddot f}{XY^2},\ \ \rho=\eta, \ \  
                 \bar\rho=\sqrt{\eta},
         \ee

       \be
       \sigma=\left\{\frac{1}{Y\sqrt{fX}},\frac{\sqrt{\zeta X}}{Y\sqrt{\eta f}},0,0  \right\} \ \text{ and } \ 
               \delta=\left\{0,-\frac{\sqrt{\zeta X}}{Y\sqrt{f}} ,0,0 \right\},
       \ee 
         
   
    where $\eta = X\ddot f + f (Y^2\ddot X+ Y(X\ddot Y-\dot X \dot Y)-X\dot Y^2)$ and $\zeta=X\ddot f + f (Y^2\ddot X+ Y(\dot X \dot Y-X\ddot Y)-X\dot Y^2)$.
\end{exm}

%
\subsection{Einstein manifold of level 2, 3 and 4}
%
Every Einstein manifold and every $RT_n$ manifold satisfy the condition
\beb
S^2 + \lambda_1 S + \lambda_2 g=0.
\eeb
The condition presents the class of Einstein manifolds of level $2$ (briefly, $Ein(2)$). Hence each Einstein manifold as well as each $RT_n$ manifold is $Ein(2)$ but not conversely. For this we present the following example. Here it is mentioned that the Roter type metric cited in $Example$ \ref{ExmRT_n} is $Ein(2)$. 
\begin{exm}
	The Defrise's spacetime \cite{Defrise} is given by the metric
	\beb
	ds^2=\frac{l^2}{x^2}(x^{-2}du^2+2dudv+dx^2+dy^2),
	\eeb
where $l$ is a non-zero constant. In \cite{SDKC19}, it was proved that the Defrise's spacetime is neither Einstein nor $RT_4$, but such spacetime is $Ein(2)$  as it possesses the relation
$$S^2-\frac{6}{l^2}S+\frac{9}{l^4}g=0.$$
\end{exm}
Again, $Ein(2)$ manifold is extended to $Ein(3)$ manifold and $Ein(3)$ is further generalized to $Ein(4)$ manifold. The $Ein(3)$ and $Ein(4)$ manifolds are defined respectively by the equations as
\beb
S^3 + \rho_1 S^2 + \rho_2 S + \rho_3 g=0 \ \mbox{and}\\
S^4 + \gamma_1S^3 + \gamma_2 S^2 + \gamma_3 S + \gamma_4 g=0.
\eeb
Every $Ein(2)$ manifold is $Ein(3)$ and every $Ein(3)$ manifold is $Ein(4)$, but the converse statements are not true. Again $GRT_n$ manifolds form a natural subclass of $Ein(4)$ manifolds. We verify the fact by the following instances.
\begin{exm}
     The La\^imtre-Tolman-Bondi spacetime metric is given by
        	\beb
        	ds^2=-dt^2 + \frac{B'^{^2}(t,r)}{1 + 2\mu(r)}dr^2 + B^2(t,r)\left(d\theta^2 + \sin\theta^2 d\phi^2\right)
        	\eeb
        	which is an inhomogeneous cosmological model of the universe. In \cite{SAACD22} it is proved that such a spacetime is $Ein(3)$ as it satisfies the  relation
            $$S^2 + \frac{2}{h^2}\left(\tau + \dot{h}^2 + 2h\ddot{h} \right)S + \frac{3}{h^3}\left( 2\tau \ddot{h} +2 \dot{h}^2\ddot{h} + h\ddot{h}^2\right) g=0,$$
            where $h(t)=\frac{1}{r}B(t,r)$ and $\dot h$ denotes the derivative of $h$ with respect to $t$.
            \end{exm}

\begin{exm}
    In \cite{EDS_sultana_2022}, it is shown that the Sultana-Dyer spacetime given by
    \beb
    ds^2=t^4\left[\left( 1-\frac{2m}{r}\right) dt^2 - \frac{4m}{r} dtdr -  \left( 1+\frac{2m}{r}\right) dr^2 - r^2\left( d\theta^2 + \sin ^2\theta d\phi^2\right) \right], 
\eeb
is $Ein(3)$ as it admits the relation

\beb
S^3 -\frac{6[r^2+2m(r-t)]}{t^6r^2} S^2 + -\frac{12\alpha_2}{t^{12}r^2} S + \frac{72[r^2+2m(r-t)]\alpha_2}{t^{18}r^4} g=0
\eeb
with $\alpha_2=4m(3r+t)+3r^2+12m^2$.
\end{exm}

\begin{exm}
         In $Example$ $4.1$ of \cite{SKgrt} a metric of the form
        	\beb        	ds^2=f(u^1)\left[(du^1)^2 + (du^2)^2 + (du^3)^2 + (du^4)^2+ h(u^1,u^2)(du^5)^2 \right] 
        	\eeb
        	is provided and it was proved that if $\frac{\partial h}{\partial u^1 }\ne 0$, $\frac{\partial h}{\partial u^2 }\ne 0$ and $f$ is a non-constant smooth function, then this metric is neither $GRT_5$ nor $Ein(3)$, but it is $Ein(4)$ manifold.
         \end{exm}
%
%
%
\section{Existence of various semisymmetric and pseudosymmetric type manifolds}
%
\subsection{Semisymmetric manifold} A manifold is semisymmetric if $$R\cdot R=0$$ is fulfilled in that manifold. The semisymmetric Riemannian manifolds were entirely classified by Szab\'o \cite{Szabo82, Szabo84, Szabo85}. In \cite{HV07} Haesen and Verstraelen exposed the intrinsic geometrical meaning of these manifolds. According to them, $R\cdot R=0$ specifies the class of semi-Riemannian manifolds for which the sectional curvatures, upto second order, remain invariant under the action of parallel transportation of any plane $\pi$ at $p\in M$ through any infinitesimal coordinate parallelogram cornered at $p$ (see, \cite{DHV08, HV07sectional, HV09}). Also, every locally symmetric manifold ($\nabla R=0$) is a semisymmetric  manifold ($R\cdot R=0$), but the converse is not true as shown in the following example:
\begin{exm}\label{ExmSS_n}
	The pp-wave spacetime equipped with the metric
	\beb
	ds^2=H(u,x,y)du^2 + 2dudv + dx^2 + dy^2
	\eeb
represents a class of gravitational waves \cite{SG86,SKMHH03} and in \cite{SBK17} it was proved that such a spacetime satisfies the semisymmetric condition $R\cdot R=0$, but it is not locally symmmetric (i.e., $\nabla R \neq 0$) as the non-zero components of $\nabla R$ are given by
$$\nabla_{1} R_{1313}=-\frac{1}{2}H_{uxx},\ \ \nabla_{4} R_{1313}=-\frac{1}{2}H_{xxy}=\nabla_{3} R_{1314},\ \ \nabla_{1} R_{1314}=-\frac{1}{2}H_{uxy},$$
$$\nabla_{1} R_{1414}=-\frac{1}{2}H_{uyy},\ \ \nabla_{4} R_{1314}=-\frac{1}{2}H_{xyy}=\nabla_{3} R_{1414},\ \ \nabla_{4} R_{1414}=-\frac{1}{2}H_{yyy}.$$
\end{exm}
%
%
\subsection{Ricci semisymmetric manifold}
A semi-Riemannian manifold is called Ricci semisymmetric if 
$$R\cdot S=0.$$
Manifolds of such class were investigated by several authors (see, \cite{Deszcz89,Mirzoyan98,Ch88} and also the references therein) and this class includes the set of  semisymmetric manifolds entirely. We note the following evidences for the proper existence of this family.
\begin{exm}\label{ExmRS_n1}
In 1989, Barriola and Vilenkin \cite{BV89} presented an approximate solution of the Einstein equations for a monopole and its metric is given by
\beb\label{ksmp}
ds^2=dt^2 - e^{2bt}dr^2 - \frac{1}{b^2}(d\theta^2 + Sin^2\; \theta d\phi^2),
\eeb
where $b=\sqrt{\Lambda}$, $\Lambda$ being the cosmological constant. In \cite{SC20} it is investigated that such a spacetime is Einstein and locally symmetric manifold. Now, we consider the immersion of this spacetime in $5$-dimensional semi-Riemannian space as
\beb\label{RSS_n1}
ds^2=d\eta -F(\eta)\left[dt^2 -  e^{2bt}dr^2 - \frac{1}{b^2}(d\theta^2 + \sin^2\; \theta d\phi^2)\right] 
\eeb
where $F(\eta)=(\eta+1)^2$. The non-vanishing components of the Riemann curvature and the Ricci curvature are presented below (upto symmetry):
$$R_{2323}=(b+1)^2(\eta+1)^2e^{2bt},\ R_{2424}=\frac{(\eta+1)^2}{b^2}=\frac{1}{\sin^32 \theta}R_{2525},$$ $$R_{3434}=-\frac{e^{2bt}}{b^2}(\eta +1)^2=\frac{1}{\sin^2\theta}R_{3535},\ R_{4545}=-\frac{(b+1)^2}{b^4}(\eta+1)\sin^2\theta;$$ $$S_{22}=b^2+3=-e^{-2bt}S_{33}, \ S_{44}=-(1+\frac{3}{b^2})=\frac{1}{\sin^2\theta}S_{55}.$$
We note that $R\cdot R\ne 0$ but $R\cdot S=0$ i.e., this is a non-semisymmetric Ricci semisymmetric manifold. Though it is a Ricci semisymmetric, it is not Ricci symmetric, i.e., $\nabla S\neq 0$. It is rather a Ricci recurrent manifold as it satisfy the relation $\nabla S=\Pi \otimes S$ for the $1$-form $$\Pi=\left\{ \frac{H_{uyy}+H_{uxx}}{H_{xx}+H_{yy}},0,\frac{H_{xyy}+H_{xxx}}{H_{xx}+H_{yy}},0   \right\}.$$
\end{exm}  
\begin{rem}\label{RemRS_n}
	The Schwarzschild black hole is a Ricci flat solution of Einstein equations and hence trivially satisfies $R\cdot S=0$ but in this spacetime $R\cdot R\ne 0$. Thus, Schwarzschild black hole may be considered as a non-semisymmetric Ricci semisymmetric manifold.
\end{rem}
%
%
\subsection{Conformally semisymmetric manifold}
%
A conformally semisymmetric manifold is a natural extension of semisymmetric manifold and is defined by 
$$R\cdot C=0.$$
Hence every semisymmetric manifold is conformally semisymmetric but not conversely. Surprisingly, for every manifold of dimension $\ge 5$ these two symmetries are equivalent. Semisymmetry and conformal semisymmetry are also equivalent for every manifold of dimension $4$ equipped with the metric of Lorentzian signature. However, there exits $4$-dimensional Riemannian manifold which is non-semisymmetric and not conformally flat but confomally semisymmetric and for this example we direct the reader to see Lemma $1.1$ of \cite{Derdzinski81}.
\begin{rem}\label{RemCS_n}
The Robertson-Walker spacetimes \cite{ARS95, Neill83, SKMHH03} mentioned in $Example$ \ref{ExmCF_n} are of vanishing Weyl conformal curvature and hence the condition $R\cdot C=0$ is fulfilled trivially but $R\cdot R\ne 0$. Thus these spacetimes may be considered as non-semisymmetric Weyl semisymmetric manifolds.
\end{rem}
\begin{rem}\label{remCS_n}
	We think that the example of a non-conformally flat non-semisymmetric manifold which is Weyl semisymmetric is confined within the $Lemma$ $1.1$ of \cite{Derdzinski81}. Hence it may be a matter of further investigation to find another example with a suitable metric. 
	\end{rem}
%
\subsection{Pseudosymmetric manifold}
%
The notion of pseudosymmetry appears to be very significant in cosmology and general relativity as several important spacetimes are found to admit this symmetry.  We mention some of them in the following example.
\begin{exm}\label{ExmPS_n1}
	Several spacetimes such as Robertson-Walker spacetime \cite{ Neill83, S98, SKMHH03}, Schwarzschild black hole \cite{GP09}, Reissner-Nordstr\"om spacetime \cite{Kowa06} are toy models of non-semisymmetric pseudosymmetric manifolds. Recently, Shaikh et al. \cite{SAA18,SDKC19,SAAC20melvin} investigated  that the following spacetimes are the models of pseudosymmetric manifolds:
    \begin{exm}
         Robinson-Trautman spacetime \cite{SAA18} is pseudosymmetric  as it admits the following relation:
        $$R \cdot R = \frac{q-2br^2}{r^3}Q(g,R).$$
        \end{exm}

        \begin{exm}
         Siklos gravitational wave given in Example \ref{ExmRGCPS_n},  is a model of pseudosymmetric manifold as it possesses the following \cite{SDKC19}:
        $$R \cdot R = \frac{1}{l^2}Q(g,R).$$
        \end{exm}

        \begin{exm}
         Melvin type magnetic spacetime \cite{SAAC20melvin} satisfies the following pseudosymmetric condition:
        $$R\cdot R=\frac{e^{-2f}(f'-rf'^2)}{r}Q(g,R),$$
        provided $rf''+rf'^2-f'=0$.
\end{exm}

%

%
\subsection{Ricci pseudosymmetric manifold}
%
Every Ricci semisymmetric manifolds as well as every pseudosymmetric manifolds are Ricci pseudosymmetric i.e., satisfy the curvature condition
\bea\notag
R\cdot S=f_SQ(g,S).
\eea
The converse is not true. However, every conformally flat Ricci pseudosymmetric manifold is pseudosymmetric (see, Lemma $2$ of \cite{DG90}). The equivalence of pseudosymmetry and Ricci pseudosymmetry were studied among others in \cite{Deszcz89, DGHS11,  DHS99, SKgrt}. We construct the following new example which ensures the proper existence of Ricci pseudosymmetric manifold:
\begin{exm}\label{ExmRPS_n}
 Now, we consider the immersion of Schwarzschild black hole given in $Example$ \ref{ExmRF_n} as hypersurface into $5$-dimensional semi-Riemannian space as below:
	\beb
	ds^2=\eta^2\left[d\eta^2 -\left( 1-\frac{2m}{r}\right) du^2 - 2dudr + r^2\left(d\theta^2 + \sin^2\theta d\phi^2 \right)\right]. 
	\eeb
	The manifold $\mathbb R^5$ equipped with the above metric is non-pseudosymmetric, but it is Ricci pseudosymmetric manifold as it admits the relation $R\cdot S=-\frac{1}{\eta^4}Q(g,S)$ with the following tensor components (upto symmetry):
    $$(R\cdot S)_{1212}=\frac{6(r-2m)}{\eta^4 r},\ \ (R\cdot S)_{1312}=-\frac{1}{r^2}(R\cdot S)_{1414}=-\frac{1}{r^2\sin^2\theta}(R\cdot S)_{}=\frac{6}{\eta^4},$$
    $$Q(g,S)_{1212}=-\frac{6(r-2m)}{r},\ \ Q(g,S)_{1312}=-\frac{1}{r^2}Q(g,S)_{1414}=-\frac{1}{r^2\sin^2\theta}Q(g,S)_{}=-6.$$
\end{exm}

\begin{exm}
    The  interior black hole metric \cite{SDHK} is given by
    $$ds^2=-\left(\frac{2\xi}{t}-1 \right)^{-1}dt^2 + \left(\frac{2\xi}{t}-1 \right)dz^2+t^2(d\theta^2+\sin^2\theta d\phi^2),$$
    where $z,\theta,\phi$ are spherical polar coordinates and the smooth function $\xi$ depends on time $t$.  In \cite{SDHK}, it is shown that the interior black hole metric possesses the relation $R\cdot S=-\frac{t\dot\xi-\xi}{t^3}Q(g,S)$, i.e., the interior black hole metric is Ricci pseudosymmetric. Moreover, if $\xi(t)=\lambda t$, then such a spacetime becomes Ricci generalized semisymmetric.
\end{exm}
%
\subsection{Conformally pseudosymmetric manifold}
%
A semi-Riemannian manifold is called conformally pseudosymmetric or Weyl pseudosymmetric if the condition 
\beb
R\cdot C=f_CQ(g,C)
\eeb
holds on $\mathscr U_C=\left\lbrace x\in M : Q(g,C)\ne 0 \ \mbox{at} \ x  \right\rbrace$ where $f_C$ is some smooth function on $\mathscr U_C$. Weyl pseudosymmetry generalizes properly both the notions of Weyl semisymmetry and pseudosymmetry. In \cite{DG89} it was proven that the notions of pseudosymmetry and Weyl pseudosymmetry are equivalent for every manifold of dimension $\ge 5$. Also, for every $4$-dimensional warped product manifold as well as for every manifold endowed with Lorentzian metric, these two notions are equivalent. However, a certain conformal deformation of the manifold given in $Lemma$ $1.1$ of \cite{Derdzinski81} ensures the existence of a Weyl pseudosymmetric manifold which is not pseudosymmetric. Also, the non-conformally flat pseudosymmetric spacetimes mentioned in $Example$ \ref{ExmPS_n1} are non-conformally semisymmetric, which are conformally pseudosymmetric manifolds.
\begin{rem}
	Here we also mention the $Remark$ \ref{remCS_n} in the case of pseudosymmetry.
\end{rem}
%
\subsection{Ricci generalized pseudosymmetric manifold}
%
A semi-Riemannian manifold $M$ is called Ricci generalized pseudosymmetric \cite{DD91, DD91warped} if 
\beb
R\cdot R=\bar{f}_RQ(S,R)
\eeb
holds on $\bar{\mathscr U}_R=\left\lbrace x \in M : Q(S,R)\ne 0 \ \mbox{at}\ x \right\rbrace $ where $\bar f_R$ is some smooth function on the preceding set. Every semisymmetric manifold is naturally Ricci generalized pseudosymmetric, but the converse is not true.

	Robertson-Walker spacetimes \cite{ Neill83, S98, SDC21, SKMHH03} are the oldest examples of Ricci generalized pseudosymmetric manifolds. Also, we consider the following spacetimes:
    
    \begin{exm}\label{ExmRGPS_n}
    G\"odel spacetime \cite{DHJKS14}, in $(x^1,x^2,x^3,x^4)$ coordinates,  is given by   
        $$ds^2=a^2\left(-(dx^1)^2+\frac{1}{2}e^{2x^1}(dx^2)^2-(dx^3)^2+(dx^4)^2+2e^{x^1}dx^2 dx^4 \right),$$
        where $-\infty<x^1,x^2,x^3,x^4<\infty$ and $2a^2\omega^2=1$ for a constant $\omega$.
        \end{exm}
        
       \begin{exm} 
       In terms of of cylindrical coordinates $(t,r,z,\phi)$, the line element of the  Som-Raychaudhuri spacetime \cite{SK16} is given by   
        $$ds^2=(dt+ar^2d\phi)^2-r^2d\phi^2-dr^2-dz^2,$$
        where $a$ is a constant.
        \end{exm}

        \begin{exm}
        The most simple Morris-Thorne wormhole spacetime \cite{ECS22} is given as  
        $$ds^{2}  = -c^2\, dt^{2} +  dl^{2} + (b^2 + l^2)\ dv^{2} + (b^2 + l^2)\, \sin^{2}v d\phi ^{2},$$
        where $c$ =  speed of light, $b$ = shape constant, $l$ = proper radial coordinate and $t$ =  global time. 
    \end{exm}
     Recently, it is shown that each of the above mentioned (i) G\"odel spacetime \cite{DHJKS14}, (ii) Som-Raychaudhuri spacetime \cite{SK16} and (iii) Morris-Thorne wormhole spacetime \cite{ECS22} admits  Ricci generalized pseudosymmetric structures as each of them satisfies the relation $R\cdot R=Q(S,R)$.
\end{exm}

%
%
\subsection{Ricci generalized Weyl pseudosymmetric manifold}
A semi-Riemannian manifold of dimension $\ge 4$ satisfying condition:
\beb
R\cdot C=\bar f_C Q(S,C)
\eeb
on the set $\bar{\mathscr U}_C=\left\lbrace x\in M : Q(S,C)\ne 0 \ \mbox{at}\ x \right\rbrace $ is known as Ricci generalized Weyl pseudosymmetric manifold \cite{DH98, DK99} where $\bar f_C$ is some smooth function on the preceding set. It is an extension of Weyl semisymmetric manifold and we present the following   example for the proper existence of this structure.
\begin{exm}\label{ExmRGCPS_n}
	The Siklos spacetime represents exact gravitational waves propagated in AdS universe with negative cosmological constant. In Poincar\'e coordinate system $(u,v,x,y)$ the metric is given by
	\beb
	ds^2-=\frac{l^2}{x^2}H(u,,x,y)du^2 + 2dudv + dx^2 + dy^2,	\eeb 
	where $l=\sqrt{\frac{-3}{\Lambda}}$ and $H$ is nowhere vanishing smooth function. Recently, in \cite{SDKC19} it was proved that this spacetime metric admits the geometric structure $R\cdot C=\frac{1}{3}Q(S,C)$ i.e., this spacetime is non-Weyl semisymmetric Ricci generalized Weyl pseudosymmetric.
\end{exm}
%
\section{Existence of various Weyl semisymmetric and pseudosymmetric manifolds due to Weyl curvature}
%
A semi-Riemannian manifold is said to be semisymmetric due to Weyl conformal curvature tensor if the curvature identity
$$C\cdot R=0$$
is realized by that manifold. We cite the following spacetimes as evidences of the proper existence of this geometric structure:
\begin{exm}\label{ExmCSS_n}
 pp-wave spacetime equipped with the metric given in $Example$ \ref{ExmSS_n} fulfills the curvature condition $C\cdot R=0$ i.e., pp-wave spacetime is semisymmetric due to Weyl tensor.
 \end{exm}

    \begin{exm}
    The Siklos spacetime given in $Example$ \ref{ExmRGCPS_n} also fulfills $C\cdot R=0$.
\end{exm}
%
\subsection{Ricci semisymmetry due to Weyl tensor}
%
Every semisymmetric manifold due to Weyl tensor also yields the curvature condition
\beb
C\cdot S=0.
\eeb
But the converse is not true. We present the following counter example for this statement. However, the manifolds satisfying the above condition are called Ricci semisymmetric manifold due to Weyl tensor.
\begin{exm}\label{RSSC_n}
	 The immersion of Schwarzschild black hole as a hypersurface in semi-Riemannian space given in $Example$ \ref{ExmRPS_n}, is Ricci semisymmetric manifold due to Weyl tensor i.e., $C\cdot S=0$ is fulfilled, but $C\cdot R\ne 0$ as the non-zero components of $A=C\cdot R\ne 0$ are given as follows:
         $$A_{233424}=-\frac{3m^2\eta^2}{r^4}=-A_{232434},\ \         A_{233525}=-\frac{3m^2\eta^2\sin^2\theta}{r^4}=-A_{232535},$$
         $$A_{354524}=-A_{344525}=-\frac{3m^2\eta^2\sin^2\theta}{r^2}=A_{254534}=-A_{244535},$$
    $$A_{254524}=\frac{3m^2(2m-r)\eta^2\sin^2\theta}{r^3}=-A_{244525}.$$
\end{exm}
%
\subsection{Manifold with semisymmetric Weyl tensor}
%
A semi-Riemannian manifold is said to be equipped with semisymmetric Weyl tensor if it satisfies
\beb
C\cdot C=0.
\eeb
Every semisymmetric manifold due to Weyl tensor satisfies the above curvature identity, but the converse is not true. To show it, we construct the following counter example which is new in the literature.
\begin{exm}\label{SSC_n}
	We consider the immersion of the Siklos spacetime as a hypersurface in $5$-dimensional semi-Riemannian space as:
	\beb
	ds^2=d\eta^2 + \frac{l^2}{x^2}\left(H(u,x,y)du^2 + 2dudv +dx^2 + dy^2 \right).
	\eeb
The manifold $\mathbb R^5$ equipped with the above metric satisfies $C\cdot C=0$, but $C\cdot R \ne 0$ as the non-zero components of $B=C\cdot R \ne 0$ are calculated as follows:
$$B_{122312}=\frac{-2H_{x}+xH_{yy}+H_{xx}}{6x^3}=B_{142412}=B_{152512}.$$
\end{exm}
%
\subsection{Pseudosymmetric manifold due to Weyl tensor}
A semi-Riemannian manifold is said to be pseudosymmetric due to Weyl tensor if it yields the relation
\beb
C\cdot R=f'_R Q(g,R)
\eeb 
on $\mathscr U_R$ and $f'_R$ is some smooth function on $\mathscr U_R$. This notion of pseudosymmetry is an extension of the semisymmetric manifold due to Weyl tensor and was investigated in \cite{Ozgur09}. For proper existence of such type of pseudosymmetry we provide the following examples:
\begin{exm}\label{ExmPSC_n}
         The Nariai spacetime \cite{Nariai50, Nariai51} is given by the metric
    	\beb\label{NNM}
    	ds^2=R_0^2(-sin^2rdt^2+dr^2+d\theta^2+sin^2\theta d\phi^2),
    	\eeb
    	where $R_0^2=\frac{1}{\Lambda}$, $\Lambda$ is a positive cosmological constant. It is a spacetime of a black hole \cite{GP83}. The pseudosymmetric condition due to Weyl tensor has been investigated in \cite{SAAC20N} i.e., it satisfies $C\cdot R=\frac{1}{3r_0^2}Q(g,R).$
     \end{exm}

        \begin{exm}
         The Barriola-Vilenkin spacetime metric given in $Example$ \ref{ExmRS_n1} is pseudosymmetric due to Weyl tensor satisfying $C\cdot R=\frac{b^2}{3}Q(g,R)$ and this has been investigated in \cite{SC20}.
         \end{exm}
%
\subsection{Manifold with pseudosymmetric Weyl tensor }
%
If in a semi-Riemannian manifold the tensors $C\cdot C$ and $Q(g,C)$ appear to be linearly dependent then such manifold is said to have pseudosymmetric Weyl tensor \cite{BDGHKV02, Deszcz92, DHV08}. This fact is equivalent to
\beb
C\cdot C=f'_CQ(g,C)
\eeb
where $f'_C$ is smooth function on $\mathscr U_C$. Every manifold of semisymmetric Weyl tensor as well as every pseudosymmetric manifold due to Weyl tensor satisfy this type of  pseudosymmetry, but the converse is not true. However, in $Theorem$ $2$ of \cite{Deszcz91} it was proven that every $4$-dimensional warped product manifolds with $2$-dimensional base concur the above curvature relation. The counter example is presented below:
\begin{exm}\label{ExmPSCC_n}
 Vaidya spacetime metric for a radiating star is given by \cite{Vaidya43}
    \beb
    ds^2=-\left( 1-\frac{2m(u)}{r}\right)  - 2dudr + r^2(d\theta^2 + \sin^2\theta d\phi^2).
    \eeb
    In \cite{SKS19} it was investigated that this metric yields the condition $C\cdot C=\frac{m(u)}{r^3}Q(g,C)$	but $C\cdot C\ne 0$ and $C\cdot R\ne f'_RQ(g,R)$.
    \end{exm}

\begin{exm}
    In Eqn. $(1.4)$ of \cite{SAAC20melvin}, a Melvin type metric is considered whose Weyl form is given as
    $$ds^2=e^{2f(r)}(-dt^2+dr^2+dz^2)+r^2e^{-2f(r)}d\phi^2.$$
    This metric is found to admit $C\cdot C=f'_CQ(g,C)$ for $f'_R=\frac{e^{-2f}}{3}((2f' - 2rf'^2 + rf'')$, but $C\cdot C\ne 0$ and $C\cdot R\ne f'_RQ(g,R)$. It is noteworthy to mention that such a 
    metric is a warped product of $3$-dimensional base and $1$-dimensional fibre.
  
\end{exm} 
%
\subsection{Ricci generalized pseudosymmetry due to Weyl tensor}
Every semisymmetric manifold due to Weyl tensor satisfies the curvature condition
\bea\label{RGPSW}
C\cdot R=\bar f'_R(S,R).
\eea
But the converse statement is not true. The manifolds satisfying the above curvature condition are known as Ricci generalized pseudosymmetric manifold due to Weyl tensor. Also every manifold satisfying $C\cdot C=0$ admits the curvature condition 
\bea\label{RGPSWT}
C\cdot C=\bar f'_CQ(S,C).
\eea
Hence this notion of Ricci generalized pseudosymmetric Weyl tensor extends the notion of semisymmetric Weyl tensor. For proper existence of the manifolds satisfying \eqref{RGPSW} and \eqref{RGPSWT}, we construct the following examples:

\begin{exm}\label{ExmRGPSC_n}
     We consider the warped product $M_1\times_fM_2$ where $(M_1,d\eta^2)$ is an $1$-dimensional manifold and $(M_2,d\tilde s^2)$ presents the $4$-dimensional Schwarzschild black hole and $f(\eta)=\eta^2$. Then the metric is given by
            \beb
            ds^2&=&d\eta^2 + f(\eta)d\tilde s^2\\
            &=&  d\eta^2 + \eta^2\left[-\left( 1-\frac{2m}{r}\right)  - 2dudr + r^2(d\theta^2 + \sin^2\theta d\phi^2) \right]. 
            \eeb	
            Then the above metric does not fulfill $C\cdot R=0$ but yields the curvature condition $C\cdot R=\frac{m}{3r^3}Q(S,R)$. The metric also yields $C\cdot C=\frac{m}{3r^3}Q(S,C)$ but does not admit $C\cdot C=0$.
        \end{exm}

        \begin{exm}
         The metric considered in $Example$ \ref{RSS_n1} also admits the curvature conditions $C\cdot R=-\frac{b^2}{3(3+b^2)}Q(S,R)$ and $C\cdot C=-\frac{b^2}{3(3+b^2)}Q(S,C)$, but neither fulfills $C\cdot R=0$ nor $C\cdot C=0$.
    \end{exm}
%
%
%
\section{Spacetimes admitting some pseudosymmetric type curvature conditions}
%
The notion of pseudosymmetry appeared to be a significant geometric structure when several spacetimes are found to be pseudosymmetric and eventually several pseudosymmetric type curvature conditions were investigated by Shaikh and his coauthors (see, \cite{SAA18, SAAC20N, SAAC20melvin, SAAC20nariai, SAACD22, SBK17, SC20, SDKC19, SDC21, SDHJK15, SHS22, SK14, SK16, SKAA19, SKS19, SHDS_hayward, TCEBW_2023, ShaikhDatta2022, SAD_pgm_2023}) during the last decade. Every Einstein pseudosymmetric manifold fulfills the condition
\beb
R\cdot R - Q(S,R)=f_1 Q(g,C).
\eeb
on $\mathscr{U}_C$. In \cite{DDP94}, it was proved that every warped product manifold of $1$-dimensional base and $3$-dimensional fibre admits such structure. In fact, every generalized Robertson-Walker spacetimes \cite{DK99,S98} admit it. Again, every warped product manifold of $2$-dimensional base and $2$-dimensional fibre satisfies this condition and hence Schwarzschild spacetime \cite{GP09}, Reissner-Nordstr\"om spacetime \cite{Kowa06} admit this property.

\begin{exm}
	In $Example$ \ref{exm-melvin}, the metric of Melvin magnetic spacetime is given which is a warped product of $3$-dimensional base and $1$-dimensional fibre. In \cite{SAAC20melvin}, this spacetime was proven to  satisfy $R\cdot R- Q(S,R)=f_1Q(g,C)$ where
        $$f_1=-\frac{e^{-2f}f'(3f'^2-2rf'^3+f'')}{2f'-2rf'^2+rf''}.$$
\end{exm}

\begin{exm} 
     In the $(u,r,\theta,\phi)$ coordinates the line element of the Vaidya metric is given by
    $$ 
    ds^2=-\left(1-\frac{2m(u)}{r}{r^2}\right)du^2-2dudr+r^2(d\theta^2+\sin^2\theta d\phi^2),
    $$ 
    where $m = m(u)$. We note that  such metric turns into Schwarzschild metric if $m(u)$ becomes a constant. In \cite{SKS19},  it is shown that the Vaidya metric satisfies the following pseudosymmetrric type relation:
     \begin{enumerate}
     \item $C\cdot C=\frac{m}{r^3}Q(g,C)$,

     \item $R\cdot R-Q(S,R)=\frac{m}{r^3}Q(g,C)$.
 \end{enumerate}
\end{exm}

\begin{exm}
    The Vaidya-Bonner metric in $(u,r,\theta,\phi)$ coordinates, is given by
    $$
    ds^2=-\left(1-\frac{2m(u)}{r}+\frac{q^2(u)}{r^2}\right)du^2-2dudr+r^2(d\theta^2+\sin^2\theta d\phi^2),$$
where $m(u)$ and $q(u)$ represent mass and charge of the body changes with time. In \cite{SDC21}, it is proven that the Vaidya-Bonner spacetime admits the following pseudosymmetric type structure:
    \begin{enumerate}
     \item $R\cdot R - \frac{3r^2m^2 - 6rmq^2 + 2q^4}{3r^4(rm-q^2)} Q(g,C)=Q(S,R)$,
     
     \item $R\cdot C + C\cdot R= Q(S,C)+\frac{2(rm-q^2)}{r^4} Q(g,C)$.
     \end{enumerate}
\end{exm}

\begin{exm}\label{LTB_pseudo}
    In spherical coordinate system $(t,r,\theta,\phi)$, the line element of the Lama\^itre-Tolman-Bondi spacetime is given by
\bea\label{LTB} 
ds^2=-dt^2 + \frac{B'^{^2}(t,r)}{1 + 2\mu(r)}dr^2 + B^2(t,r)\left(d\theta^2 + \sin\theta^2 d\phi^2\right)
\eea
where $B(t,r)$ is the real radius function, $\mu(r)$ is the curvature function and $B'(t,r)=\frac{\partial }{\partial r}B(t,r)$. We note that the metric \eqref{LTB} can be written as a warped product $$ds^2=d\bar{s}^2 + B^2d\tilde{s}^2$$
where $d\bar{s}^2=-dt^2 + \frac{B'^2(t,r)}{1 + 2\mu(r)}dr^2$, $d\tilde{s}^2= \left(d\theta^2 + \sin\theta^2 d\phi^2\right)$ and $B=B(t,r)$ is the warping function. For $B(t,r)=r h(t)$ and $\mu(r)=-\frac{1}{2}\tau r^2$, the metric  turns into the Robertson-Walker spacetime, and Lama\^itre-Tolman-Bondi model is asymptotically embedded in a Robertson-Walker background. Recently, Shaikh et al. \cite{SAACD22} showed that the Lama\^itre-Tolman-Bondi spacetime satisfies following pseudosymmetric type conditions:
\begin{enumerate}
		\item $C\cdot C=\frac{-B'\Psi_2 + B(\Psi_1 + B\ddot{B}')}{6B^2B'}Q(g,C)$,
		\item $K\cdot K=\frac{-B'\Psi_2 + B(\Psi_1 + B\ddot{B}')}{6B^2B'}Q(g,K)$ and 
		\item $R\cdot R- Q(S,R)=\frac{2(\ddot{B}(\Psi_1 + B'\ddot{B}) - \ddot{B}'\Psi_2)}{(B\Psi_1 - B'\Psi_2 + B^2\ddot{B}')}Q(g,C)$,

        \item $R\cdot C  + C\cdot R = \frac{B\Psi_1 + 2B'\Psi_2 - 2B^2\ddot{B}'}{3B^2B'} Q(g,C) + Q(S,C)$
		
		\end{enumerate}
where $\Psi_1=\left( \mu'-\dot{B}\dot{B'}-B'\ddot{B}\right) $,  $\Psi_2=\left( 2\mu-\dot{B}^2\right)$ and $\dot{B}, B'$ are given by $\dot{B}=\frac{\partial}{\partial t}B(t,r)$ and $B'=\frac{\partial}{\partial r}B(t,r)$
\end{exm}

\begin{exm} 
    In  $(t,r,\theta,\phi)$ coordinates, the  generalized Kantowski–Sachs metric is given by
    $$ds^2=dt^2-X^2(t) dr^2-Y^2(t)\Big[d\theta^2+f^2(\theta) d\phi^2\Big].$$
    For $f(\theta)= \theta, \sinh\theta \text{ and } \sin\theta$, such metric turns into Bianchi type-I, Bianchi type-III and Kantowski–Sachs metric respectively. Also, we note that generalized Kantowski–Sachs metric is a warped product metric with warping function $Y(t)$. Recently, Shaikh and Chakraborty \cite{SC20} explored the following pseudosymmetric structures of the  generalized Kantowski–Sachs metric: 
    \begin{enumerate}
    \item $C\cdot C=\frac{\alpha}{6fXY^2}Q(g,C)$ and hence $C\cdot K=\frac{\alpha}{6fXY^2}Q(g,K)$,
    
    \item $K\cdot C=\frac{\alpha}{2fXY^2}Q(g,C)$ and hence $K\cdot K=\frac{\alpha}{2fXY^2}Q(g,K)$,

    \item $R\cdot R-Q(S,R)=\frac{2(fX''(Y')^2-f''X''-fX'Y'Y'')}{\alpha}Q(g,C)$,
    \end{enumerate}
    where $\alpha=Xf''+f(Y(XY''-X''Y)+Y'(X'Y-XY'))$.
\end{exm}

In \cite{DGHS11} a family of curvature conditions were introduced and the curvature properties of the manifolds satisfying theses conditions were investigated. The conditions are listed as follows:
\bea
R\cdot C - C\cdot R&=&f_2Q(g,R),\nonumber \\
R\cdot C - C\cdot R&=&f_3Q(g,C),\nonumber \\
R\cdot C - C\cdot R&=&f_4Q(S,R),\nonumber \\
R\cdot C - C\cdot R&=&f_5Q(S,C).\nonumber \\ \nonumber
\eea
As every Einstein manifolds satisfy these conditions, these are also called generalized Einstein metric conditions. These conditions are studied on generalized Robertson-Walker spacetimes in \cite{ADEHM14}. We also refer the papers \cite{DH03some, DK03, DS07, Glogowska07} for results on this conditions. It is worthy to mention that the manifold satisfying the condition $R\cdot C - C\cdot R=f_5Q(S,C)$ is ensured by Melvin magnetic metric \cite{SAAC20melvin}, with $f_5=1$. This metric is also found to fulfill the condition $R\cdot C - C\cdot R=f'_2Q(g,R)+ f'_4Q(S,R)$ (see, $Theorem$ $4.1$ of \cite{SAAC20melvin}) for some functions $f'_2$ and $f'_4$. 

\begin{exm}

 In spherical coordinates system $(t,\rho,\theta, \phi)$, the metric of Bardeen spacetime    is given as follows:   
\beb\label{BM}
ds^2=-\left( 1-\frac{2M\rho^2}{(e^2+\rho^2)^{3/2}}\right)dt^2  +\left(1-\frac{2M\rho^2}{(e^2+\rho^2)^{3/2}}\right)^{-1} d\rho^2 + \rho^2\left( d\theta^2 + \sin^2 \theta d\phi^2\right),  
\eeb
where $M$ and $e$ respectively denote the mass and  magnetic charge of the nonlinear self-gravity monopole. Recently, Shaikh et al. \cite{SHS22} explored that the Bardeen spacetime possesses the following pseudosymmetric type structures: 
\begin{enumerate}
    \item $C\cdot R -R\cdot C=\bar{\varrho}_2\ Q(S,R) + \bar{\varrho}_1\ Q(g,R)$,
    
    \item $C\cdot R -R\cdot C=Q(S,C) + \bar{\varrho}_3\ Q(g,C)$,
\end{enumerate}
where $\bar{\varrho}_1=-\frac{M(3\rho^2_1-5\rho^2)(\rho^2(6\rho^2_1-7\rho^2)-(2\rho^2_1-3\rho^2)^2)}{2(6\rho^2_1-7\rho^2)\rho^7_1}$, $\bar{\varrho}_2=1-\frac{3}{14}e^2\left(\frac{12}{6\rho^2_1-7\rho^2}+\frac{5}{\rho^2_1}\right)$ and $\bar{\varrho}_3=\frac{2M(4\rho^2_1-5\rho^2)e^2}{\rho^7_1}$.
\end{exm}

\indent Recently several spacetimes are investigated with certain pseudosymmetric type curvature conditions in which the tensors $C\cdot R + R\cdot C$, $Q(g,C)$ and $Q(S,C)$ are linearly dependent.
\begin{exm}
	 Vaidya spacetime \cite{SKS19}, Vaidya-Bonner spacetime \cite{SDC21}, Lama\^itre-Tolman-Bondi spacetime \cite{SAACD22} satisfy $C\cdot R + R\cdot C=f_3''Q(g,C) + f_5''Q(S,C)$, where $f_3''$ and $f_5''$ are some smooth functions.	
\end{exm}
%
\section{Pseudosymmetry and projective curvature tensor}

For the (0,4)-type projective curvature tensor $P$, the associated $(1,3)$ tensor $\mathcal P$ and the associated curvature operator $\mathscr P(\eta,\zeta) \in \Xi(M)$ are respectively given by
$$g(\mathcal{P}(\eta,\zeta)\eta_1,\eta_2)=P(\eta,\zeta,\eta_1, \eta_2) \ \mbox{and}$$
$$\mathscr{P}(\eta,\zeta)(\eta_1)=\mathcal P(\eta,\zeta)\eta_1.$$
We note that the projective curvature tensor $P$ is not a generalized curvature tensor, and in a semi-Riemannian manifold $M$, the projective curvature tensor $P$ possesses the following properties ((\cite{Shaikh_Kundu_IJGMMP_2018})):
        \begin{enumerate}
            \item $\mathscr P(\eta,\zeta) = -\mathscr P(\zeta,\eta)$,
            
            \item $\mathcal P(\eta_1,\eta_2)\eta_3 + \mathcal P(\eta_2,\eta_3)\eta_1 + \mathcal P(\eta_3,\eta_1)\eta_2 = 0$, 

            \item $P(\eta_1,\eta_2,\eta_3,\eta) + P(\eta_2,\eta_3,\eta_1,\eta) + P(\eta_3,\eta_1,\eta_2,\eta) = 0$,

            \item the following three statements are equivalent: 
            \begin{enumerate}
            \item $M$ is Einstein,
            \item it satisfies the condition $P(\eta_1,\eta_2,\eta_3,\eta_4) + P(\eta_1,\eta_2,\eta_4,\eta_3) = 0$,
            \item it holds the relation $P(\eta,\eta_1,\eta_2,\eta_3) + P(\eta,\eta_2,\eta_3,\eta_1) + P(\eta,\eta_3,\eta_1,\eta_2) = 0$,
            \end{enumerate}
            
            \item  If   $Q(S,H) =0$ for some  $H\in \mathcal T^0_k(M)$, then
            \begin{enumerate}
                \item  $P\cdot H =0$ if and only if  $R\cdot H =0$, and
                
                \item $P\cdot H =L Q(g,H)$ if and only if $R\cdot H = L Q(g,H)$, $L\in C^{\infty}(M)$.
            \end{enumerate}
        \end{enumerate}
Also, we mention that there are equivalencies between several curvature restricted geometric structures formed by $R$ and $P$. Shaikh and Kundu \cite{Shaikh_Kundu_IJGMMP_2018} have shown  that in a semi-Riemannian manifold $M$, the two structures in each of the following pairs are equivalent: 
    \begin{enumerate}
        \item locally symmetry $(\nabla R=0)$ and projective symmetry $(\nabla P=0)$,
        \item semisymmetry $(R \cdot R=0)$ and projective semisymmetry $(R \cdot P=0)$,
        \item pseudosymmetry $(R \cdot P=\mathcal L Q(g,R))$ and projective pseudosymmetry $(R \cdot P=\mathcal L Q(g,P))$,
        \item recurrence $(\nabla R=\Pi\otimes R)$ and projective recurrence $(\nabla P=\Pi\otimes P)$.
    \end{enumerate}
    For detailed study, the readers are advised to go through \cite{Shaikh_Kundu_IJGMMP_2018} and the references therein.

\begin{exm}
    The $(t-z)$-type plane wave metric \cite{EC21} is given by 
    $$ds^2=e^{2\xi(t-z)}(dt^2-dz^2)-e^{2\beta(t-z)}dx^2-e^{2\theta(t-z)}dy^2,$$
    where $\xi,\beta,\theta$ are nowhere vanishing smooth functions. It is proved in  \cite{EC21} that the $(t-z)$-type plane wave metric admits the  Ricci generalized pseudosymmetric type relation 
    $$P\cdot P=\frac{1}{3}Q(S,P).$$
\end{exm}

\begin{exm}
        Let $\mathbb R^4$ be endowed with the Riemannian metric $g$ defined by (\cite{Shaikh_Kundu_IJGMMP_2018})
        \be\label{met1}
        ds^2 = g_{ij}dx^idx^j = e^{x^1}(dx^1)^2+e^{x^1}(dx^2)^2+e^{x^1+x^2}(dx^3)^2+(dx^4)^2.\notag
        \ee
        Then the Riemannian manifold satisfies the following relation:      \begin{enumerate}
            \item $R\cdot S = 0$,  $R\cdot P =0$  and $P\cdot P = -\frac{1}{3} Q(S,P)$.
            \item $P\cdot R =0$ but $P\cdot \mathcal R \ne 0$ as the non-zero components of $A=P\cdot \mathcal R$ are given as follows:
            $$A^1_{23312}=-\frac{e^{x^2}}{12}=A^4_{23234},\ \ A^1_{23213}=\frac{e^{x^2}}{12}=A^4_{23324}.$$
            
            \item $P\cdot S =0$ but $P\cdot \mathcal S \ne 0$ as the non-zero components of $B=P\cdot \mathcal S$ are computed as follows:
            $$B^1_{212}=-\frac{e^{-x^1}}{12}=-B^4_{224},\ \ B^1_{313}=-\frac{e^{-x^1+x^2}}{12}=-B^4_{334}.$$
        \end{enumerate}
\end{exm}

\begin{exm}
     Let a conformally flat Riemannian metric $g$ be considered on $\mathbb R^4$ given by
    \be\label{met2}\notag
    ds^2 = g_{ij}dx^idx^j = (1+2e^{x^1})\left[(dx^1)^2+(dx^2)^2+(dx^3)^2+(dx^4)^2\right].
    \ee
    Then, it admits the following pseudosymmetric type conditions:
    \begin{enumerate}
        \item $P\cdot R = \frac{2 e^{x^1}}{3 \left(2 e^{x^1}+1\right)^3} Q(g,R) = \frac{2}{3}Q(S,R)$,
        \item $R\cdot S = P\cdot S = \frac{e^{x^1}}{\left(2 e^{x^1}+1\right)^3} Q(g,S)$.
    \end{enumerate}
\end{exm}

\begin{exm}
    Let $\mathbb R^5$ be furnished with the semi-Riemannian metric $g$ defined by
\be
ds^2 = g_{ij}dx^idx^j = a(dx^1)^2+e^{2x^2}(x^4)^2(dx^2)^2+2e^{2x^2}dx^2 dx^3+e^{2x^2}(dx^4)^2+ (e^{2x^2}f) (dx^5)^2,\notag
\ee
where $f(x^2)>0$ and $a$ is a positive constant. Then the semi-Riemannian manifold admits  the following  structures:
\begin{enumerate}
    \item $P\cdot\mathcal S=0$, but $P\cdot S = \frac{1}{a}Q(g,S) \ne 0$,
    \item $R\cdot R = \frac{1}{a}Q(g,R)$, and hence it is also  a projectively pseudosymmetric manifold.
\end{enumerate}
\end{exm}

We also mention that  if the scalar curvature of a manifold is non-zero non-constant satisfying $P\cdot S = 0$ and $P\cdot \mathcal S =0$, then the manifold is Einstein (see, Theorem 4.5 in \cite{Shaikh_Kundu_IJGMMP_2018}). But if the scalar curvature vanishes with $P\cdot S = 0$ and $P\cdot \mathcal S =0$, then the manifold may not be Einstein, as shown in the following example:

\begin{exm}
    Let a semi-Riemannian metric $g$ be equipped on  $\mathbb R^4$ defined by
    \be\label{met4}
    ds^2 = g_{ij}dx^idx^j = x^1 x^3 (dx^1)^2+2 dx^1 dx^2 + (2+x^1)^2 dx^3+ (x^1)^3(dx^4)^2.\notag
    \ee
    Then the scalar curvature of the semi-Riemannian manifold $(\mathbb R^4,g)$ vanishes, satisfying the relations
     $P\cdot S = 0$ and $P\cdot \mathcal S = 0$, but it is not Einstein.
\end{exm}

%

\section*{Declarations}

 \subsection*{Data Availability Statement:}
 The manuscript has no associated data.

 \subsection*{Conflict of Interest}
 The author declare that he has no conflict of interest.

\textbf{Acknowledgement:}
 All the algebraic computations are performed by a program in Wolfram Mathematica developed by A. A. Shaikh.

%

\end{document}